\documentclass[11pt]{article}

\usepackage{theorem,amssymb,amsmath}

\usepackage{graphicx}
\usepackage{caption}

\usepackage{hyperref}

\advance \topmargin by -\headheight
\advance \topmargin by -\headsep
\evensidemargin \oddsidemargin
\author{
Jean-Paul Allouche \\
CNRS, IMJ-PRG, Sorbonne \\
4 Place Jussieu \\
F-75252 Paris Cedex 05, France \\
\url{https.www.imj-prg.fr/~jean-paul.allouche}
\and
Doron Zeilberger \\
Department of Mathematics \\
Rutgers University (New Brunswick) \\
Hill Center-Busch Campus, 110 Frelinghuysen Rd. \\
Piscataway, NJ 08854-8019, USA \\
\url{https://sites.math.rutgers.edu/~zeilberg/}
\and
with an appendix by Shalosh B. Ekhad \\
{\tt ShaloshBEkhad@gmail.com}}

\title{Human and automated approaches for finite trigonometric sums}

\date{ }

\def \proof{\bigbreak\noindent{\it Proof.\ \ }}

\def \endpf{{\ \ $\Box$ \medbreak}}

\newtheorem{theorem}{Theorem}
\newtheorem{lemma}{Lemma}
\newtheorem{proposition}{Proposition}

\theorembodyfont{\rm}

\newtheorem{remark}{Remark}

\begin{document}

\maketitle

\begin{quotation}\begin{em}
\flushright To the memory of Vladimir Shevelev (Mar 09, 1945 -- May 03, 2018)
\par
\bigskip
\end{em}\end{quotation}

\begin{abstract}
We show that identities involving trigonometric sums recently proved by Harshitha, Vasuki and
Yathirajsharma, using Ramanujan's theory of theta functions, were either already in the literature
or can be proved easily by adapting results that can be found in the literature. Also we prove
two conjectures given in that paper. After mentioning many other works dealing with identities 
for various trigonometric sums, we end this paper by describing an automated approach for proving 
such trigonometric identities.

\medskip

\end{abstract}

\medskip

Identities between quantities that are not clearly and immediately related are always fascinating.
In particular when their proofs use subtle or complicated arguments (think of the Basel Problem,
and the solution by Euler, $\sum 1/n^2 = \pi^2/6$). But it also happens that, after a ``complicated''
proof is found, people try to find a ``direct'' or an ``elementary'' proof: in some cases this is not (not 
yet?) possible (think of the Fermat-Wiles theorem). Trigonometric identities, {\it a priori} unexpected,
can also be found in the literature. We will just cite two such identities, one due to Gauss
(see, e.g., \cite{gauss}, where some signs seem to be misprinted):
$$
\cos\left(\frac{\pi}{17}\right) =
\frac{1 - \sqrt{17} + \sqrt{34-2\sqrt{17}} + 2 \sqrt{17 + 3\sqrt{17} + \sqrt{34-2\sqrt{17}}
+2\sqrt{34+2\sqrt{17}}}}{16}
$$
the other is a consequence of Ramanujan theta function identities (take the first two equalities 
of Lemma~3.8 in \cite{be-zh}, subtract the second one from the square of the first one, or see 
\cite[Corollary~32.1,~p.\ 184]{berndt}; also see Eq. (1.31) of Corollary~7 in \cite{liu}, and 
Eq.~(1.11) in \cite{ha-va-ya}):
$$
\frac{\sin^2\left(\frac{3\pi}{7}\right)}{\sin\left(\frac{2\pi}{7}\right)}
- \frac{\sin^2\left(\frac{2\pi}{7}\right)}{\sin\left(\frac{\pi}{7}\right)}
+ \frac{\sin^2\left(\frac{\pi}{7}\right)}{\sin\left(\frac{3\pi}{7}\right)} = 0.
$$

\medskip

Much more modestly we stumbled upon the recent paper \cite{ha-va-ya} published in the 
Ramanujan journal, where the authors give an interesting way to obtain a closed form
for several trigonometric sums by using Ramanujan's theory of theta functions. In particular
Theorem~1.1 in that paper gives six formulas about which the authors write that these six identities
{\it seem to be new}. Of course it is very interesting to obtain these equalities through theta functions, 
but we were interested to search whether these formulas previously appeared in the literature.
Having found a vast set of papers on the subject of identities involving finite trigonometric sums,
we were lucky enough to discover either references in which some of these identities had already 
been proved, or known results from which the other identities can be easily deduced. Here we 
describe or adapt these results. We also confirm two conjectures given by the authors of 
\cite{ha-va-ya} at the end of their paper.

\medskip

Finally we show how to give automated proofs for trigonometric identities of this kind.

\section{The six identities of Harshitha, Vasuki and Yathirajsharma}

First we recall the six trigonometric identities given in \cite[Theorem~1.1]{ha-va-ya}.

\bigskip

\begin{itemize}

\item[ ]{1.}
{\it If $k$ is an odd natural number, then}
\begin{equation}\label{Eq05}
\sum_{j=1}^{\frac{k-1}{2}} (-1)^{j-1} \sin\left(\frac{(2j-1)\pi}{2k}\right) = \frac{(-1)^{\frac{k-3}{2}}}{2}
\ \ \ \text{\it (here one should suppose that $k \geq 3$)}
\end{equation}
and
\begin{equation}\label{Eq06}
\sum_{j=1}^{\frac{k-1}{2}} (-1)^{j-1} \csc\left(\frac{(2j-1)\pi}{2k}\right) = 
\frac{k + (-1)^{  \frac{k+1}{2}}}{2}\cdot
\end{equation}

\item[ ]{2.}
{\it If $n$ is an odd natural number and $j=2p$ is an even positive integer such that
$\gcd(j,n) = 1$, then}
\begin{equation}\label{Eq07}
\sum_{k=0}^{\frac{n-3}{2}}
\frac{\sin\left(\frac{(j+1)(2k+1)\pi}{2n}\right)\sin\left(\frac{(j-1)(2k+1)\pi}{2n}\right)} 
{\sin^2\left(\frac{(2k+1)\pi}{2n}\right)\sin^2\left(\frac{p(2k+1)\pi}{n}\right)}
= \frac{n^2 - 1}{3}\cdot
\end{equation}

\item[ ]{3.}
{\it If $n$ is an even number and $j \equiv 2 \pmod 4$, such that $\gcd(\frac{j}{2}, \frac{n}{2}) = 1$, 
then}
\begin{equation}\label{Eq08}
\sum_{k=0}^{\frac{n}{2}-1}
\frac{\sin\left(\frac{(j+1)(2k+1)\pi}{2n}\right)\sin\left(\frac{(j-1)(2k+1)\pi}{2n}\right)} 
{\sin^2\left(\frac{(2k+1)\pi}{2n}\right)\sin^2\left(\frac{p(2k+1)\pi}{n}\right)}
= \frac{n^2}{4}\cdot
\end{equation}

\item[ ]{4.}
{\it If $n$ is an even number and $j \equiv 2 \pmod 4$, such that $\gcd(\frac{j}{2}, \frac{n}{2}) = 1$, 
then}
\begin{equation}\label{Eq09}
\sum_{k=0}^{\frac{n}{2}-1}
\frac{\sin\left(\frac{(j+1)k\pi}{n}\right)\sin\left(\frac{(j-1)k\pi}{n}\right)} 
{\sin^2\left(\frac{k\pi}{n}\right)\sin^2\left(\frac{jk\pi}{n}\right)}
= \frac{n^2 - 4}{12}\cdot
\end{equation}

\item[ ]{5.}
{\it If $n$ is an odd natural number and $j$ is a positive integer such that $\gcd(j,n) = 1$, then}
\begin{equation}\label{Eq10}
\sum_{k=0}^{\frac{n - 3}{2}}
\frac{\sin\left(\frac{(j+1)(2k+1)\pi}{n}\right)\sin\left(\frac{(j-1)(2k+1)\pi}{n}\right)} 
{\sin^2\left(\frac{(2k+1)\pi}{n}\right)\sin^2\left(\frac{j(2k+1)\pi}{n}\right)}
= 0.
\end{equation}

\end{itemize}

\section{Some related trigonometric identities found in the literature}

Here we give three propositions that state results we found in the literature.

\begin{proposition}[{\rm \cite{pa-re}, last displayed formula}]\label{cf01}
$$
\sum_{j=1}^n (-1)^{j+1} \sin\left(\frac{(2j-1)\pi}{4n+2}\right) = \frac{(-1)^{n+1}}{2}\cdot
$$
\end{proposition}

\begin{remark}\label{same-kind}
Several sums of the same kind appear in the literature, where $\sin$ can be replaced 
with $\cos$ or a power of $\cos$ or $\sin$, for example a companion formula of the formula 
above is given in \cite{pa-re}, namely, for $n \geq 1$,
$$
\sum_{j=1}^n (-1)^{j+1} \cos\left(\frac{j\pi}{2n+1}\right) = \frac{1}{2}\cdot
$$
Also the relations
$$
\sum_{k=1}^n (-1)^{k+1} \cos^2\left(\frac{k\pi}{2n+2}\right) = \frac{1}{2} \ \ \ \text{and} \ \ \
\sum_{
k=1}^n (-1)^{k+1} \sin^2\left(\frac{k\pi}{2n+2}\right) = \frac{(-1)^{n+1}}{2}
$$
can be found in \cite{weller}. The same formulas can be found in \cite[Theorem~3.4, p.~218]{an-ha}.

\medskip

\noindent
Note that there is a straightfoward proof of the summation in Proposition~\ref{cf01} as well as of the 
three summations in this remark, by using the relation (see, e.g., \cite[Example~7, p.~405]{chrystal}, 
or use the complex exponential form of cos and sin)
$$
(\#) \ \ \ \ \ \cos(\alpha) + \cos(\alpha+\beta) + \dots + \cos(\alpha+(K-1)\beta) 
= \frac{\sin\left(\frac{1}{2}\beta K\right)}{\sin\left(\frac{1}{2}\beta\right)} 
\cos\left(\alpha + \frac{1}{2} \beta (K-1)\right)
$$
and replacing the possible occurrences of $\cos^2(x)$ and $\sin^2 x$ with $(1 \pm \cos 2x)/2$.
\end{remark}

\begin{proposition}[{\rm \cite{hassan}, Relation~3.16 on p.~817}]\label{cf02}
$$
\sum_{j=0}^{n-1} (-1)^j \csc\left(\frac{(2j+1)\pi}{4n+2}\right) = 
n + \frac{1 - (-1)^n}{2}\cdot
$$
\end{proposition}

Now we give the last proposition of this section that we discovered in a 1908 book by Bromwich. 

\begin{proposition}[{\rm \cite{bromwich}, p.\ 183; also see p.\ 187}]\label{bro}
$$
\left\{\begin{array}{llclcl}
&\text{\rm (I)} \ \ \ \ \ &\displaystyle\sum_{k=1}^{\frac{n-1}{2}} \frac{1}{\sin^2\left(\frac{k\pi}{n}\right)} 
&=& \dfrac{n^2-1}{6} &\ \ \text{if $n$ is odd;} \\
\\
&\text{\rm (II)} \ \ \ \ \ &\displaystyle\sum_{k=1}^{\frac{n}{2}-1} \frac{1}{\sin^2\left(\frac{k\pi}{n}\right)} 
&=& \dfrac{n^2-4}{6} &\ \ \text{if $n$ is even;} \\
\\
&\text{\rm (III)} \ \ \ \ \ &\displaystyle\sum_{k=0}^{\frac{n-3}{2}} \frac{1}{\sin^2\left(\frac{(2k+1)\pi}{2n}\right)} 
&=& \dfrac{n^2-1}{2} &\ \ \text{if $n$ is odd;} \\
\\
&\text{\rm (IV)} \ \ \ \ \ &\displaystyle\sum_{k=0}^{\frac{n}{2}-1} \frac{1}{\sin^2\left(\frac{(2k+1)\pi}{2n}\right)} 
&=& \dfrac{n^2}{2} &\ \ \text{if $n$ is even.} \\
\end{array}
\right.
$$
\end{proposition}

\begin{remark}\label{fish}
One can find in \cite[Identity~(23)]{fisher} the identity 
$$
\sum_{k=1}^{n-1} \frac{1}{\sin^2\left(\frac{k\pi}{n}\right)} = \frac{n^2-1}{3}\cdot
$$
The case where $n$ is odd is an exercise in Bromwich's book (\cite[Ex.\ 4, p.\ 188]{bromwich}).
Also see \cite{gr-ru-lo} where it is stated that such identities were conjectured while computing 
low-temperature series for a $Z_n$-symmetric Hamiltonian in statistical mechanics \cite{he-la}; 
(also see \cite{he} and  the references therein for much more on this and similar identities). 
Actually this identity can be easily deduced from (I) and (II) in Proposition \ref{bro} above
(cut the sum into $\displaystyle\sum_{k=1}^{(n-1)/2} + \displaystyle\sum_{(n+1)/2}^{n-1}$ if 
$n$ is odd, and into $\displaystyle\sum_{k=1}^{n/2} + \displaystyle\sum_{(n/2)+1}^{n-1}$ if $n$
is even, and make the change of index $\ell = n - k$ in each of the second sums).

\end{remark}

\begin{remark}\label{bru}

One can find in \cite{bruch} (also see the other references given there) the identity
$$
\sum_{k=0}^{n-1} \frac{1}{\sin^2\left(\frac{(2k+1)\pi}{2n}\right)} = n^2.
$$
Note that this equality is given in \cite[Remark~2]{ha-va-ya} under the slightly disguised form
$$
\sum_{k=1}^{n} \frac{1}{\sin^2\left(\frac{(2k+1)\pi}{2n}\right)} = n^2
$$
as a corollary of the identity
$$
\sum_{k=1}^{n-1} \cot^2\left(\frac{k\pi}{n}\right) = \frac{(n-1)(n-2)}{3} \ \ \text{for any $n > 0$.}
$$
The authors of \cite{ha-va-ya} cite a similar derivation in \cite[Corollary~2.4]{be-ye}. It is 
interesting to note that the ``disguised'' equality above is already given (actually only for $n$ even) 
under a slightly different disguise by Riesz in 1914 in \cite{riesz1, riesz2} where it is used to prove 
a theorem of Bernstein (about Riesz' results also see \cite[Volume II, p.~10,~Eq.~3.11]{zygmund} 
and \cite{chu}). Actually the identity at the beginning of this Remark~\ref{bru} can be easily 
deduced from (III) and (IV) in Proposition \ref{bro} above, as was done in Remark~\ref{fish}.
\end{remark}

\begin{remark}
Note that the identities (I, II, III, IV) in Proposition~\ref{bro} are, up to notation, respectively the
identities (3.40, 3.28, 3.14, 3.4) in the paper \cite{hassan}, where they are also generalized.
\end{remark}

\medskip

We end this section with an easy lemma partly deduced from what precedes.

\begin{lemma}\label{easy}
The following equalities hold.

\begin{itemize}

\item[ ]{\text{(i)}}
Let $n$ be an odd integer and $p$ an integer such that $\gcd(p,n) = 1$. Then
$$
\sum_{k=0}^{\frac{n-3}{2}} \frac{1}{\sin^2\left(\frac{p(2k+1)\pi}{n}\right)}
= \sum_{\ell = 1}^{n-1} \frac{1}{\sin^2\left(\frac{p \ell \pi}{n}\right)} 
- \frac{1}{2} \sum_{r = 1}^{n-1} \frac{1}{\sin^2\left(\frac{2 p r \pi}{n}\right)}
= \frac{n^2-1}{6}\cdot
$$

\item[ ]{\text{(ii)}}
Let $n$ be an even integer and $p$ an odd integer, such that $\gcd(p,n/2) = 1$. Then
$$
\sum_{k=1}^{\frac{n}{2}-1}\frac{1}{\sin^2\left(\frac{p(2k+1)\pi}{n}\right)}
= \sum_{\ell = 1}^{n-1} \frac{1}{\sin^2\left(\frac{p \ell \pi}{n}\right)} -
\sum_{r = 1}^{\frac{n}{2}-1} \frac{1} {\sin^2\left(\frac{2 p r \pi}{n}\right)}
= \frac{n^2}{4}\cdot
$$

\end{itemize}

\end{lemma}

\proof
$(i)$ For $n$ odd and $\gcd(p,n) = 1$ we have
$$
\begin{array}{lll}
\displaystyle\sum_{k=0}^{\frac{n-3}{2}} \frac{1}{\sin^2\left(\frac{p(2k+1)\pi}{n}\right)} 
&=& \displaystyle\sum_{\substack{1 \leq \ell \leq n-2 \\ \ell \ \text{odd}}} 
\frac{1}{\sin^2\left(\frac{p \ell \pi}{n}\right)} \\
&=& \displaystyle\sum_{1 \leq \ell \leq n-1} \frac{1}{\sin^2\left(\frac{p \ell \pi}{n}\right)}
- \sum_{\substack{1 \leq \ell \leq n-1 \\ \ell \ \text{even}}} \frac{1}{\sin^2\left(\frac{p \ell \pi}{n}\right)} \\
&=&  \displaystyle\sum_{1 \leq \ell \leq n-1} \frac{1}{\sin^2\left(\frac{p \ell \pi}{n}\right)}
- \sum_{1 \leq r \leq \frac{n-1}{2}} \frac{1}{\sin^2\left(\frac{2p r \pi}{n}\right)}\cdot
\end{array}
$$
But 
$$
\sum_{1 \leq r \leq \frac{n-1}{2}} \frac{1}{\sin^2\left(\frac{2p r \pi}{n}\right)}
= \sum_{\frac{n+1}{2} \leq s \leq n-1}  \frac{1}{\sin^2\left(\frac{2p s \pi}{n}\right)} \ \
\text{(change of index $s = n - r$)}.
$$
Hence
$$
\sum_{1 \leq r \leq \frac{n-1}{2}} \frac{1}{\sin^2\left(\frac{2p r \pi}{n}\right)} =
\frac{1}{2}\left(\sum_{1 \leq r \leq \frac{n-1}{2}} \frac{1}{\sin^2\left(\frac{2p r \pi}{n}\right)} 
+  \sum_{\frac{n+1}{2} \leq s \leq n-1}  \frac{1}{\sin^2\left(\frac{2p s \pi}{n}\right)}\right)
= \frac{1}{2} \sum_{1 \leq r \leq n-1} \frac{1}{\sin^2\left(\frac{2p r \pi}{n}\right)}\cdot
$$
Thus, finally
$$
\sum_{k=0}^{\frac{n-3}{2}} \frac{1}{\sin^2\left(\frac{p(2k+1)\pi}{n}\right)}
= \sum_{1 \leq \ell \leq n-1} \frac{1}{\sin^2\left(\frac{p \ell \pi}{n}\right)}
- \frac{1}{2} \sum_{1 \leq r \leq n-1} \frac{1}{\sin^2\left(\frac{2p r \pi}{n}\right)}\cdot
$$

Now we claim that the first sum on the right hand side only depends on the value of $p \ell$ modulo 
$n$ and that $p \ell$ takes exactly once all the values modulo $n$ except the value $0$ since 
$\gcd(p, n) = 1$, thus this sum is equal to the same sum where $p$ is replaced with $1$. 
Similarly the second sum only depends on the value of $2pr$ modulo $n$ and $2pr$ takes exactly 
once all the values modulo $n$ except the value $0$ since $\gcd(2p,n)=1$, thus the sum is equal to 
the same sum where $2p$ is replaced with $1$. Thus
$$
\sum_{k=0}^{\frac{n-3}{2}} \frac{1}{\sin^2\left(\frac{p(2k+1)\pi}{n}\right)} =
 \sum_{1 \leq \ell \leq n-1} \frac{1}{\sin^2\left(\frac{\ell \pi}{n}\right)}
- \frac{1}{2} \sum_{1 \leq r \leq n-1} \frac{1}{\sin^2\left(\frac{r \pi}{n}\right)}
= \frac{1}{2} \sum_{1 \leq r \leq n-1} \frac{1}{\sin^2\left(\frac{r \pi}{n}\right)}\cdot
$$
This last sum is equal to $\frac{n^2-1}{6}$ from Remark~\ref{fish}, where it was
indicated that the first identity there is a consequence of Proposition~\ref{bro}~(I).

\bigskip

$(ii)$ For $n$ even, $p$ odd and $\gcd(p,n/2) = 1$ we have
$$
\begin{array}{lll}
\displaystyle\sum_{k=0}^{\frac{n}{2}-1} \frac{1}{\sin^2\left(\frac{p(2k+1)\pi}{n}\right)} 
&=& \displaystyle\sum_{\substack{1 \leq \ell \leq n-1 \\ \ell \ \text{odd}}} 
\frac{1}{\sin^2\left(\frac{p \ell \pi}{n}\right)} \\
&=& \displaystyle\sum_{1 \leq \ell \leq n-1} \frac{1}{\sin^2\left(\frac{p \ell \pi}{n}\right)}
- \sum_{\substack{1 \leq \ell \leq n-1 \\ \ell \ \text{even}}} \frac{1}{\sin^2\left(\frac{p \ell \pi}{n}\right)} \\
&=&  \displaystyle\sum_{1 \leq \ell \leq n-1} \frac{1}{\sin^2\left(\frac{p \ell \pi}{n}\right)}
- \sum_{1 \leq r \leq \frac{n}{2}-1} \frac{1}{\sin^2\left(\frac{2p r \pi}{n}\right)}\cdot
\end{array}
$$
As seen above the first sum on the right hand side does not depend on $p$. It is equal to
$\frac{n^2-1}{3}$. Similarly, defining $m := n/2$, the second sum can be written
$$
\sum_{1 \leq r \leq m-1} \frac{1}{\sin^2\left(\frac{p r \pi}{m}\right)}
$$
where $\gcd(p,m)=1$. It is thus equal to $\frac{m^2-1}{3}$. Hence, finally
 $$
\sum_{1 \leq k \leq \frac{n}{2}-1} \frac{1}{\sin^2\left(\frac{p (2k+1) \pi}{n}\right)} = \frac{n^2}{4}\cdot
$$

\begin{remark}\label{in-passing}
In the course of the proofs above, we proved, for $\gcd(p,n) = 1$, the identity
$$
\sum_{1 \leq \ell \leq n-1} \frac{1}{\sin^2\left(\frac{p \ell \pi}{n}\right)} = \frac{n^2-1}{3}
$$
(compare with Remark~\ref{fish}).
\end{remark}

\section{Identities (\ref{Eq05}) to (\ref{Eq10}) revisited}

We begin this section with a straightforward trigonometric formula, namely
$$
(\#) \ \ \ \ \  \frac{\sin x \sin y}{\sin^2(\frac{x+y}{2}) \sin^2(\frac{x-y}{2})} = 
\frac{1}{\sin^2 (\frac{x-y}{2})} - \frac{1}{\sin^2 (\frac{x+y}{2})} 
$$
(write: $2 \sin x \sin y = \cos(x-y) - \cos(x+y)$ and use $\cos z = 1 - 2 \sin^2(z/2)$).

\medskip

Now we revisit Identities (\ref{Eq05}) to (\ref{Eq10}), to show that they were either
already in the literature or easy consequences of identities in the literature.

\begin{itemize}

\item Identity~(\ref{Eq05})

See Proposition~\ref{cf01}.

\item Identity~(\ref{Eq06})

See Proposition~\ref{cf02}.

\item Identity~(\ref{Eq07})

We have for $n$ odd, $j=2p$ and $\gcd(j,n)=1$, successively using $(\#)$, 
Proposition~\ref{bro}~(III), and Lemma~\ref{easy}~(ii),
$$
\begin{array}{lll}
\displaystyle\sum_{k=0}^{\frac{n-3}{2}}
\frac{\sin\left(\frac{(j+1)(2k+1)\pi}{2n}\right)\sin\left(\frac{(j-1)(2k+1)\pi}{2n}\right)} 
{\sin^2\left(\frac{(2k+1)\pi}{2n}\right)\sin^2\left(\frac{p(2k+1)\pi}{n}\right)} 
&=& \displaystyle\sum_{k=0}^{\frac{n-3}{2}} \frac{1}{\sin^2\left(\frac{(2k+1)\pi}{2n}\right)}
- \sum_{k=0}^{\frac{n-3}{2}} \frac{1}{\sin^2\left(\frac{p(2k+1)\pi}{n}\right)} \\
&=& \dfrac{n^2-1}{2} - \dfrac{n^2-1}{6} = \dfrac{n^2-1}{3}\cdot
\end{array}
$$
\item Identity (\ref{Eq08})

We have for $n$ even, $j=2p$, $j \equiv 2 \pmod 4$, and $\gcd(\frac{j}{2}, \frac{n}{2}) = 1$
successively using $(\#)$, Proposition~\ref{bro}~(IV), and Lemma~\ref{easy}~$(i)$
$$
\begin{array}{lll}
\displaystyle\sum_{k=1}^{\frac{n}{2}-1}
\frac{\sin\left(\frac{(j+1)(2k+1)\pi}{2n}\right)\sin\left(\frac{(j-1)(2k+1)\pi}{2n}\right)} 
{\sin^2\left(\frac{(2k+1)\pi}{2n}\right)\sin^2\left(\frac{p(2k+1)\pi}{n}\right)}
&=& \displaystyle\sum_{k=1}^{\frac{n}{2}-1} \frac{1}{\sin^2\left(\frac{(2k+1)\pi}{2n}\right)}
- \sum_{k=1}^{\frac{n}{2}-1} \frac{1}{\sin^2\left(\frac{p(2k+1)\pi}{n}\right)} \\
&=& \dfrac{n^2}{2} - \dfrac{n^2}{4} = \dfrac{n^2}{4}\cdot
\end{array}
$$

\item Identity (\ref{Eq09})

We have for $n$ even, $j=2p$, $j \equiv 2 \pmod 4$, and $\gcd(j/2,n/2)=1$, successively using
$(\#)$, Proposition~\ref{bro}~(II), and the identity in Remark~\ref{in-passing},
$$
\begin{array}{lll}
\displaystyle\sum_{k=1}^{\frac{n}{2}-1}
\frac{\sin\left(\frac{(j+1)k\pi}{n}\right)\sin\left(\frac{(j-1)k\pi}{n}\right)} 
{\sin^2\left(\frac{k\pi}{n}\right)\sin^2\left(\frac{jk\pi}{n}\right)}
&=& \displaystyle\sum_{k=1}^{\frac{n}{2}-1} \frac{1}{\sin^2\left(\frac{k\pi}{n}\right)}
- \displaystyle\sum_{k=1}^{\frac{n}{2}-1} \frac{1}{\sin^2\left(\frac{jk\pi}{n}\right)} \\
&=& \displaystyle\sum_{k=1}^{\frac{n}{2}-1} \frac{1}{\sin^2\left(\frac{k\pi}{n}\right)}
- \displaystyle\sum_{k=1}^{\frac{n}{2}-1} \frac{1}{\sin^2\left(\frac{(j/2)k\pi}{n/2}\right)} \\
&=& \dfrac{n^2-4}{6} - \dfrac{(\frac{n}{2})^2-1}{3} = \dfrac{n^2-4}{12}\cdot
\end{array}
$$

\item Identity (\ref{Eq10})

We have for $n$ odd and $j$ such that $\gcd(j,n)=1$, successively using $(\#)$, 
and Lemma~\ref{easy}$(i)$,
$$
\begin{array}{lll}
\displaystyle\sum_{k=0}^{\frac{n - 3}{2}}
\frac{\sin\left(\frac{(j+1)(2k+1)\pi}{n}\right)\sin\left(\frac{(j-1)(2k+1)\pi}{n}\right)} 
{\sin^2\left(\frac{(2k+1)\pi}{n}\right)\sin^2\left(\frac{j(2k+1)\pi}{n}\right)}
&=& \displaystyle\sum_{k=0}^{\frac{n - 3}{2}}\frac{1}{\sin^2\left(\frac{(2k+1)\pi}{n}\right)}
-  \displaystyle\sum_{k=0}^{\frac{n - 3}{2}}\frac{1}{\sin^2\left(\frac{j(2k+1)\pi}{n}\right)} \\
&=& 0.
\end{array}
$$

\end{itemize}

\section{Confirming two conjectures of Harshitha, Vasuki, Yathirajsharma}

At the end of their paper \cite{ha-va-ya}, the authors give two conjectures about which they
wrote: {\it We feel that these two can also be tackled with Ramanujan's theory although we 
were unable to do so.} First we recall these conjectures.

\medskip

\noindent
{\bf Conjecture \cite{ha-va-ya}} \  {\it If $k$ is any positive integer, then
$$
\lim_{k \to \infty} \left\{\sum_{j=0}^{2k-1} (-1)^j \sin^2\left(\frac{(2j+1)\pi}{8k+2}\right)\right\}
= - \frac{1}{2}
$$
and
$$
\lim_{k \to \infty}  \left\{\sum_{j=0}^{2k} (-1)^j \sin^2\left(\frac{(2j+1)\pi}{8k+6}\right)\right\} 
= \frac{1}{2}\cdot
$$
}

Here we give a simple proposition (but without Ramanujan's theory of theta functions) which 
easily implies both conjectures above. 

\begin{proposition}\label{simple}
The following equalities hold.
$$
\sum_{j=0}^{2k-1} (-1)^j \sin^2\left(\frac{(2j+1)\pi}{8k+2}\right) = 
\frac{- \sin^2\left(\dfrac{2k\pi}{4k+1}\right)}{2\cos\left(\dfrac{\pi}{4k+1}\right)}
$$
and
$$
\sum_{j=0}^{2k} (-1)^j \sin^2\left(\frac{(2j+1)\pi}{8k+6}\right) = 
\frac{1}{2} - \frac{\cos^2\left(\dfrac{(2k+1)\pi}{4k+3}\right)}{2\cos\left(\dfrac{\pi}{4k+3}\right)}\cdot
$$
As a consequence the conjectures above hold.
\end{proposition}

\proof The proof is left to the reader who can use the relations
$\sin^2(x) = \frac{1}{2}(1 - \cos(2x))$, $(-1)^t \cos y = \cos(t\pi + y)$, and, 
as previously, Identity~$(\#)$  at the end of Remark~\ref{same-kind}. 
\endpf

\section{More on identities for trigonometric sums}

Actually the literature on identities for finite trigonometric sums is huge, from ancient references 
like \cite{ramus, eisenstein, stern} to papers published this year, but here we have only scratched 
the surface of this vast literature: transitively searching through the references that we give, as well 
as tracking other papers (that we have skipped because they often deal with weighted sums) with, 
e.g.,  the keyword ``cotangent sums'', yield a lot of many other interesting items. We will only cite
\cite{ej-le} (and the references therein).

\medskip

Such trigonometric sums occur in a large variety of domains, and are addressed with a
large variety of methods: compiling all the pertinent papers in the domain would lead at least
to a book (look for example at the long list of references given in \cite{be-ye}). Various methods
have been used, some of them in papers cited above: {\it interpolation formulas} 
\cite{riesz1, riesz2, an-ha, hassan} (to which we can add, e.g., \cite{by-sm, an-as}...), 
{\it Ramanujan's theta function} \cite{be-zh, ha-va-ya, vi-sh-ha}, {\it calculus of residues} (see, 
among several other papers, the papers by Cvijovi\'c or Cvijovi\'c et al. cited in \cite{cvijovic}; also
see \cite{gra-pro}), {\it expansions in partial fractions} (see, e.g., \cite{chu-mar, wang-zheng, chu}), 
{\it discrete Fourier analysis} (see \cite{be-ha}), etc. 

\medskip

The applications go from {\it number theory} (see, e.g., \cite{by-sm}, but also
papers dealing with Dedekind sums, generalizations and associated reciprocity laws (see, e.g.,
\cite{zagier, fukuhara, be-za}...) to {\it physics} (see, in particular, 
\cite{dowker1, dowker2, cv-sr1, cv-sr2} and the references therein; also see \cite{holcombe}), from 
{\it enumerative combinatorics}  (e.g., number of closed walks on a path as in \cite{daf-gla-kow1}) 
to {\it binomial identities} (see, e.g. \cite{ramus, merca2012-ro, merca2012, merca2014, holcombe}), 
and to {\it topology} (see, e.g., \cite{lawson}), etc. 

\bigskip

To end this section, we would like to cite two more results. The first one is a relation between 
Dirichlet series with trigonometric coefficients and finite trigonometric sums given, among other 
results, in a paper of J. Franke \cite{franke}.

\bigskip

\noindent{\bf Theorem} (\cite{franke}, particular case) \ The following relation holds
$$
\sum_{\stackrel{\scriptsize n > 0}{\scriptsize N \nmid n}} \frac{\cot^2(\frac{n \pi}{N})}{n^2} =
\frac{(N-1)(N-2)(N^2+3N+2)\pi^2}{90 N^2}\cdot
$$  
The second one, that we find particularly nice, is a result of V. Shevelev \cite{she} and V. Shevelev 
and P. J. C. Moses \cite{she-mos}, linking tangent power sums and sums of digits of integers in 
integer bases.
 
\bigskip

\noindent {\bf Theorem} (\cite{she, she-mos}) \ 
{\it Let $n$ be an odd integer $\geq 3$. Let $s_{n-1}(r)$ be the sum of digits of the integer $r$
in base $n-1$. Define $S_n$ by
$$
S_n(x) = \sum_{\stackrel{\scriptstyle 0 \leq r  < x}{\scriptstyle r \equiv 0 \bmod{n}}} (-1)^{s_{n-1}(r)}.
$$
Then if $p$ is a positive integer, one has
$$
S_n((n-1)^{2p}) = \frac{2}{n} \
\sum_{k=1}^\frac{n-1}{2} \tan^{2p} \left(\frac{\pi k}{n}\right)
\sim \ \frac{2}{n} (n-1)^{2p \lambda_n} \ \ \text{as $p$ tends to infinity,}
$$
where $\lambda_n = \dfrac{\log \cot (\frac{\pi}{2n})}{\log(n-1)}\cdot$
}

\begin{remark}
This last result shows why $\sum_{k=1}^\frac{n-1}{2} \tan^{2p} \left(\frac{\pi k}{n}\right)$ is an integer
multiple of $n$. Furthermore, for $n = 3$, the reader can see that $\lambda_3 = \frac{\log 3}{\log 4}$ 
and guess that there is a nice relation to the Moser-Newman statement which observes then
quantifies that, among the integers multiple of $3$, there are ``more'' evil than odious numbers 
(recall that a number is said to be {\it evil} if the sum of its binary digits is even and {\it odious} 
otherwise): everything is explained in \cite{she, she-mos} and the references therein. 
\end{remark}

\section{Proving and reproving}

As indicated in \cite{be-ye} {\it This history} [the history of evaluations and reciprocity theorems 
for trigonometric sums] {\it is, not surprisingly, sporadic, and consequently authors often publish 
results without being aware that their theorems had previously been published elsewhere.} 
Given that this was written 20 years ago, it should not be a surprise that more recent papers
give examples strengthening this remark. We give two more examples: first, compare, say, the abstract 
of \cite{go-za} with Section~3 of \cite{hassan} (or even with the abstract of \cite{gau-bru}); second, 
we will briefly look at some of the results obtained in \cite[Theorem~2]{vi-sh-ha} via identities for 
the Ramanujan theta function, and show that they are either in the literature or they can be readily 
obtained from results already in the literature. First we state an easy result.

\begin{theorem}\label{extra-ident}
We have, for $n$ odd, the following identities.
$$
\sum_{k=1}^{\frac{n-1}{2}} (-1)^{k+1} \frac{\sin \frac{2k\pi}{n}}{\sin\frac{k\pi}{n}} = 1
$$
and
$$
\sum_{k=1}^{\frac{n-1}{2}} (-1)^{k+1}\frac{\sin\frac{k\pi}{n}}{\sin\frac{2k\pi}{n}} = 
(-1)^{\frac{n+1}{2}} \left(\frac{n-1}{4}\right) + \frac{\chi_{o}\left(\frac{n-1}{2}\right)}{2}
$$
where $\chi_{o}(m) = 1$ if $m$ is odd and $\chi_{o}(m) = 0$ if $m$ is even.
\end{theorem}

\proof
The first equality is equivalent to: for odd $n$ one has
$$
2 \sum_{k=1}^{\frac{n-1}{2}} (-1)^{k+1} \cos \frac{k\pi}{n} = 1,
$$
which is given Remark~\ref{same-kind} above.

The second equality is an easy consequence of the identity (for odd $n$)
$$
\sum_{k=1}^{\frac{n-1}{2}} \frac{(-1)^{k+1}}{\cos\frac{k\pi}{n}} = 
(-1)^{\frac{n+1}{2}} \left(\frac{n-1}{2}\right) + \chi_{o}\left(\frac{n-1}{2}\right)
$$
which can be found, e.g., in \cite{hassan} (Identity~(3.43)).

\begin{remark}
The first identity in Theorem~\ref{extra-ident} above gives the identities (1.1), (1.3), and (1.8) 
of \cite[Theorem~2]{vi-sh-ha} (take $n = 7, 13, 17$ and use, if needed, that $\sin(\pi-x) = \sin x$). 
Note that, as indicated in \cite{vi-sh-ha}, their identity (1.1) is already in \cite{liu}. The second 
equality in Theorem~\ref{extra-ident} above gives the identity (1.2) of \cite[Theorem~2]{vi-sh-ha} 
(take $n = 7$), which, as indicated in \cite{vi-sh-ha}, is already in \cite{liu}.

\medskip

Other identities of \cite[Theorem~2]{vi-sh-ha} can be proved in an elementary way: for example
identities (1.4) and (1.9), which have the form $\sum \pm (\sin a \ \sin b)/(\sin c \ \sin d)$. We only 
show how to prove easily (1.4). Let $L$ be the lefthand term of Identity~(1.4), i.e.,
$$
L := \frac{\sin\frac{4\pi}{13}}{\sin\frac{2\pi}{13}}\frac{\sin\frac{6\pi}{13}}{\sin\frac{3\pi}{13}}
- \frac{\sin\frac{2\pi}{13}}{\sin\frac{\pi}{13}}\frac{\sin\frac{3\pi}{13}}{\sin\frac{5\pi}{13}}
- \frac{\sin\frac{5\pi}{13}}{\sin\frac{4\pi}{13}}\frac{\sin\frac{\pi}{13}}{\sin\frac{6\pi}{13}}
$$
which can also be written, using the relations $\sin x = \sin(\pi - x)$, $\sin(2x) = 2 \sin x \cos x$,
$2 \cos x \cos y = \cos(x+y) + \cos(x-y)$ and $\cos(\pi - x) = - \cos x$,
$$
\begin{array}{lll}
L &:=& \displaystyle\frac{\sin\frac{4\pi}{13}}{\sin\frac{2\pi}{13}}\frac{\sin\frac{6\pi}{13}}{\sin\frac{3\pi}{13}}
- \frac{\sin\frac{2\pi}{13}}{\sin\frac{\pi}{13}}\frac{\sin\frac{10\pi}{13}}{\sin\frac{5\pi}{13}}
- \frac{\sin\frac{8\pi}{13}}{\sin\frac{4\pi}{13}}\frac{\sin\frac{12\pi}{13}}{\sin\frac{6\pi}{13}} \\
&  & \\
&=& \displaystyle 4 \cos \frac{2\pi}{13} \cos \frac{3\pi}{13} - 4 \cos \frac{\pi}{13} \cos \frac{5\pi}{13} 
- 4 \cos \frac{4\pi}{13} \cos \frac{6\pi}{13} \\
&  & \\
&=& \displaystyle 2 \cos \frac{5\pi}{13}  + 2\cos \frac{\pi}{13} - 2 \cos \frac{6\pi}{13} - 2 \cos \frac{4\pi}{13}
- 2 \cos \frac{10\pi}{13} - 2 \cos \frac{2\pi}{13} \\
&  & \\
&=& \displaystyle 2 \sum_{k=1}^6 (-1)^{k+1} \cos \frac{k\pi}{13} = 1 \ \ \text{(see again
Remark~\ref{same-kind} above)} \\
\end{array}
$$
Actually one can note that, going backwards and grouping differently the terms in the last 
sum just above (possibly using $\cos(\pi-x) = - \cos x$ and $\sin(\pi -x) = \sin x$), yields other 
identities, for example:
$$
\begin{array}{lll}
1 &=& \displaystyle 2 \sum_{k=1}^6 (-1)^{k+1} \cos \frac{k\pi}{13}  \\
   &  & \\
   &=& \displaystyle 2 \left(\cos \frac{\pi}{13} + \cos \frac{3\pi}{13}\right) 
          - 2 \left(\cos \frac{6\pi}{13} + \cos \frac{2\pi}{13}\right)
           + 2 \left(\cos \frac{5\pi}{13} - \cos \frac{4\pi}{13}\right) \\
   & & \\       
   &=& \displaystyle 2 \left(\cos \frac{\pi}{13} + \cos \frac{3\pi}{13}\right) - 
          2 \left(\cos \frac{6\pi}{13} + \cos \frac{2\pi}{13}\right)
           - 2 \left(\cos \frac{8\pi}{13} + \cos \frac{4\pi}{13}\right) \\
    & & \\         
    &=& \displaystyle 4 \left(\cos \frac{2\pi}{13} \cos \frac{\pi}{13}\right) 
            -4 \left(\cos \frac{4\pi}{13} \cos \frac{2\pi}{13}\right)
            -4 \left(\cos \frac{6\pi}{13} \cos \frac{2\pi}{13}\right) \\  
     & & \\
     &=&  \displaystyle\frac{\sin\frac{4\pi}{13}}{\sin\frac{2\pi}{13}}\frac{\sin\frac{2\pi}{13}}{\sin\frac{\pi}{13}}
- \frac{\sin\frac{8\pi}{13}}{\sin\frac{4\pi}{13}}\frac{\sin\frac{4\pi}{13}}{\sin\frac{2\pi}{13}}
- \frac{\sin\frac{12\pi}{13}}{\sin\frac{6\pi}{13}}\frac{\sin\frac{4\pi}{13}}{\sin\frac{2\pi}{13}} \\       
\end{array}
$$
which can be ``disguised'' by replacing $\sin x$ with $\sin(\pi - x)$, or here by simplifying and
using $\sin(\pi - x) = \sin x$ again, thus obtaining
$$
\frac{\sin \frac{4\pi}{13}}{\sin \frac{\pi}{13}} - \frac{\sin\frac{8\pi}{13}}{\sin\frac{2\pi}{13}}
- \frac{\sin\frac{\pi}{13}}{\sin\frac{6\pi}{13}}\frac{\sin\frac{4\pi}{13}}{\sin\frac{2\pi}{13}} = 1, \ \text{etc.}
$$
\end{remark}

\section{Automated proofs of the above and similar identities}\label{auto}

Given the type of identities described above, and the re-discoveries of some of them, it seems
desirable to have automated proofs: two examples are \cite{daf-gla, daf-gla-kow2} 
where the authors allude to the use of {\tt Mathematica}.
But it also seems desirable to have not only automated proofs but also automated ``discoveries''
(and proofs).

Two hints seem to be in favor of possible automated discoveries-proofs (or ``shaloshable 
discovery-proofs'', see \cite{zeilberger}; also see \cite{PWZ}) of the above and similar trigonometric 
identities. The (very vague) first hint is the occurrence of a hypergeometric series in 
\cite[Theorem~2.1]{ha-va-ya}. The second hint is the following identity that can be found, 
e.g., in \cite[p.\ 184]{bromwich}

\begin{equation}\label{prodsin}
\sin(n\theta) = 2^{n-1}\prod_{r=0}^{n-1} \sin\left(\theta+\frac{r\pi}{n}\right)\cdot
\end{equation}

\noindent
This identity can be used to prove a bunch of identities. For example, taking its logarithmic derivative,
then differentiating, yields
$$
\frac{n^2}{\sin^2(n \theta)} - \frac{1}{\sin^2(\theta)} 
= \sum_{r=1}^{n-1} \frac{1}{\sin^2\left(\theta + \frac{r\pi}{n}\right)}\cdot
$$
Now, letting $\theta$ tend to $0$, one obtains the identity of Remark~\ref{fish}.
But one way to obtain Identity~(\ref{prodsin}) is to use Cebyshev polynomials.
Thus one could imagine that manipulating Cebyshev polynomials and the like
with a computer algebra system, conveniently taught, would permit to discover/prove
trigonometric identities ``automatically''.

\subsection{Automated proofs of specific finite trigonometric identities}

Given a multivariable polynomial $P$, any identity of the form 
$$
P\left(\sin\left(\frac{\phantom{2}\!\!\!\pi}{n}\right), \, \sin\left(\frac{2 \pi}{n}\right), \dots, \, \sin\left(\frac{(n-1)\pi}{n}\right)\right) = 0 
\eqno(1)
$$
for a {\it specific}, positive integer $n$, is routinely provable by {\tt Maple}, and (probably) any other 
computer algebra system. In {\tt Maple} the command is {\tt simplify}.

\medskip

\noindent
For example the original `Morrie's law' \cite{morrie}, that Richard Feynman never forgot (also see 
\cite{lo-be-ze}),
$$
\cos(20^{\circ}) \, \cos(40^{\circ}) \, \cos(80^{\circ}) \, = \, \frac{1}{8}
$$
can be automatically done (exactly!) by {\tt Maple}. Just type
$$
{\tt simplify(cos(Pi/9)*cos(2*Pi/9)*cos(4*Pi/9));}
$$
and you will get right away {\tt 1/8}.

\medskip

\noindent
Nowadays you do not have to be a Gauss to (rigorously!) prove Gauss' identity at the beginning of this 
article, just enter, in {\tt Maple}
\begin{multline*}
{\tt
simplify((1-sqrt(17)+sqrt(34-2*sqrt(17))+ 2*sqrt(17+3*}\backslash \\
{\tt sqrt(17)+sqrt(34-2*sqrt(17))+2*sqrt(34+2*sqrt(17))))/16-cos(Pi/17));
}
\end{multline*}
and in one nano-second you would get {\tt 0}. 

\medskip

\noindent
To take another example, to prove the second identity on page 1, type
$$
{\tt sin(3/7*Pi)**2/sin(2/7*Pi) - sin(2/7*Pi)**2/sin(1/7*Pi) + sin(1/7*Pi)**2/sin(3/7*Pi):} 
$$
$$
{\tt simplify(\%);} \phantom{sin(3/7*Pi)**2/sin(2/7*Pi) - sin(2/7*Pi)**2/sin(1/7*Pi) + sin(1/7*Pi)**2}
$$
and you would immediately get  {\tt 0}.

\medskip

\noindent
The way {\tt Maple} does it is to use Euler's formula
$$
\sin x = \frac{e^{ix}-e^{-ix}}{2i}\cdot
$$
Then in $(1)$ everything is a polynomial in $z:=e^{i\pi/n}$, getting, upon expansion, a polynomial 
in $z$ and making sure that it is divisible by $z^n+1$. So numeric computations turn into routine, 
symbolic,  `high-school algebra' calculations that computer algebra systems excel at.

\medskip

\noindent
Nine such identities are proved, by `advanced' methods in \cite{vi-sh-ha}. See the output file
$$
\text{\url{http://www.math.rutgers.edu/~zeilberg/tokhniot/oTrigSums14.txt}}
$$
for automatic, elementary ({\it high-school algebra}) proofs of all of them.

\subsection{Automated proofs of infinite families of finite trigonometric sums}

More interesting are ``infinite families'' of finite trigonometric sums. We will describe how to 
automatically derive explicit polynomial expressions to the following eight families.

\medskip

Below $n$ and $k$ are  arbitrary positive integers.

\medskip

{\bf Type Top}
$$
\sum_{j=1}^{n} \sin^{2k} \left(\frac{j \pi}{2n+1}\right).
$$

{\bf Type Ton}
$$
\sum_{j=1}^{n} \csc^{2k} \left(\frac{j \pi}{2n+1}\right).
$$

{\bf Type Tep}
$$
\sum_{j=1}^{n} \sin^{2k} \left(\frac{(2j-1) \pi}{4n}\right).
$$

{\bf Type Ten}
$$
\sum_{j=1}^{n} \csc^{2k} \left(\frac{(2j-1) \pi}{4n}\right).
$$

{\bf Type Uop}
$$
\sum_{j=1}^{n-1} \sin^{2k} \left(\frac{j \pi}{2n}\right).
$$

{\bf Type Uon}
$$
\sum_{j=1}^{n-1} \csc^{2k} \left(\frac{j \pi}{2n}\right).
$$

{\bf Type Uep}
$$
\sum_{j=1}^{n} \sin^{2k} \left(\frac{(2j-1) \pi}{4n+2}\right).
$$

{\bf Type Uen}
$$
\sum_{j=1}^{n} \csc^{2k} \left(\frac{(2j-1) \pi}{4n+2}\right).
$$

It turns out that to get closed-form expressions for general $n$ and $k$ for the four types {\bf Top, Tep, 
Uop, Uep}, (i.e., for sums of {\bf positive} even powers of the sines), one does not need
computers. These are easy human exercises. Let us just illustrate it with type {\bf Top}.

\medskip

Define $z:=e^{\frac{i \pi}{2n+1}}$. Note that $z^{2n+1}=-1$. Then  
$\sin(\frac{j \pi}{2n+1})=(z^j-z^{-j})/(2i)$. We have:
$$
\sum_{j=1}^{n} \sin^{2k} \left(\frac{j \pi}{2n+1}\right) \, = \,
\sum_{j=1}^{n}  \left ( (z^j -z^{-j})/(2i) \right )^{2k} \, = \,
\frac{(-1)^k}{4^k} \sum_{j=1}^{n} (z^j -z^{-j})^{2k}.
$$
By the binomial theorem, this equals
$$
\frac{(-1)^k}{4^k} \sum_{j=1}^{n} \sum_{r=0}^{2k} (-1)^{r} \frac{(2k)!}{r!(2k-r)!} z^{-jr+j(2k-r)} \, = \,
\frac{(-1)^k}{4^k} \sum_{j=1}^{n} \sum_{r=0}^{2k} (-1)^{r} \frac{(2k)!}{r!(2k-r)!} z^{2j(k-r)}
$$
$$
\begin{aligned}
&= \, \frac{(-1)^k}{4^k} \cdot n \cdot (-1)^{k} \frac{(2k)!}{k!^2} +
\frac{(-1)^k}{4^k} \sum_{j=1}^{n} \sum_{r=0}^{k-1} (-1)^{r} \frac{(2k)!}{r!(2k-r)!} (z^{2j(k-r)}+z^{-2j(k-r)}) \\
&=\, \frac{n}{4^k} \frac{(2k)!}{k!^2} +
\frac{(-1)^k}{4^k}  \sum_{r=0}^{k-1} (-1)^{r}  \frac{(2k)!}{r!(2k-r)!} \sum_{j=1}^{n} (z^{2j(k-r)}+z^{-2j(k-r)}). 
\end{aligned}
$$
Note that, by summing the geometric series and using $z^{2(2n+1)}=1$, we have:
$$
 \sum_{j=1}^{n} (z^{2j(k-r)}+z^{-2j(k-r)})=
 \sum_{j=-n}^{n} (z^{2(k-r)})^j -1 \, = \, 0-1 \, = \, -1.
$$
Going back to the above, we have that our desired sum equals
$$
=\frac{n}{4^k} \frac{(2k)!}{k!^2} +\frac{(-1)^{k+1}}{4^k}  \sum_{r=0}^{k-1} (-1)^{r}  \frac{(2k)!}{r!(2k-r)!}\cdot
$$
But (you prove it!)
$$
 \sum_{r=0}^{k-1} (-1)^{r}  \frac{(2k)!}{r!(2k-r)!}  = \frac{1}{2} (-1)^{k+1} \frac{(2k)!}{k!^2} \cdot
$$
We have just proved, {\it sans ordinateurs}, the general identity:

\medskip

\noindent{\bf Theorem Top}

$$
\sum_{j=1}^{n} \sin^{2k} \left(\frac{j\, \pi}{2n+1}\right) \, = \,
\frac{1}{4^k} \frac{(2k)!}{k!^2} \left(n+ \frac{1}{2} \right).
$$

\medskip

Similar things can be done with the other three positive families, and they are left to the reader.

\medskip

\noindent{\bf Theorem Tep}
$$
\sum_{j=1}^{n} \sin^{2k} \left(\frac{(2j-1) \pi}{4n}\right) \, = \, 
\frac{n}{2^{2k-1}} \, \frac{(2k-1)!}{(k-1)!k!} \cdot
$$

\medskip

\noindent{\bf Theorem Uop}
$$
\sum_{j=1}^{n-1} \sin^{2k} \left(\frac{j \pi}{2n}\right) 
\, = \, \frac{n}{4^k} \frac{(2k)!}{k!^2} - \frac{1}{2}  \cdot
$$

\medskip

\noindent{\bf Theorem Uep}
$$
\sum_{j=1}^{n} \sin^{2k} \left(\frac{(2j-1) \pi}{4n+2}\right) 
\, = \,
\frac{n}{4^k} \frac{(2k)!}{k!^2} +
\frac{1}{4^k} \frac{(2k-1)!}{k!(k-1)!} -\frac{1}{2} \cdot
$$

\medskip

Note that in the above four theorems we have {\bf closed-form} expressions  for {\it general} 
$n$ {\bf and}  {\it general} $k$.

\medskip

Much more interesting is to derive expressions for the sum of the even powers of the 
{\bf cosecants} (i.e., reciprocals of the sines). 

It turns out that for the case {\bf Uon}  Chu and Marini 
\cite{chu-mar} (bottom of Page~126) found an elegant, but very complicated, explicit formula that yields, for each specific $k$,
a polynomial expression, in $n$, of degree $2k$.  Our approach, to be described shortly, while not `explicit' (it uses recursion), does the same, and seems far more efficient.
Chu and Marini also gave an elegant generating function, whose output agrees with ours.

Everything is implemented in the {\tt Maple} package {\tt TrigSums.txt} available, along with numerous input and output files from
$$
\text{\url{https://sites.math.rutgers.edu/~zeilberg/mamarim/mamarimhtml/trig.html}}.
$$

\begin{remark}
The results for $k \leq 60$ obtained through this package are given at 
$$
\text{\url{https://sites.math.rutgers.edu/~zeilberg/tokhniot/oTrigSums10.txt}}
$$
(also see {\bf Type Uon} in the Appendix) while, using the generating function of \cite{chu-mar}, 
the results are given at
$$
\text{\url{https://sites.math.rutgers.edu/~zeilberg/tokhniot/oTrigSums10ChuMarini.txt}}
$$
which takes about the same time: note that, using the explicit formula in \cite{chu-mar} is
much slower.
\end{remark}

\medskip

\noindent
We will illustrate our approach with the family of type {\bf Ton}. 
Recall that the {\bf Chebyshev polynomial of the first kind}, $T_n(x)$ may be defined by
$$
T_n(\cos t) = \cos(n t).
$$
It is well-known and easy to see \cite{chebyshev}, that

$$
T_{n}(x) \, = \, \frac{n}{2} \sum_{k=0}^{\lfloor \frac{n}{2} \rfloor} \frac{(-1)^k 2^{n-2k} (n-k-1)!}{k! (n-2k)!} x^{n-2k}.
$$
Hence
$$
\frac{T_{2n+1}(x)}{x} \, = \, (2n+1) \sum_{k=0}^{n} \frac{(-1)^k 4^{n-k} (2n-k)!}{k! (2n+1-2k)!} (x^2)^{n-k}.
$$
Define the degree $n$ polynomial
$$
E_n(x):=\frac{T_{2n+1}(\sqrt{x})}{\sqrt{x}}.
$$
Then
$$
E_n(x) \, = \, (2n+1) \sum_{k=0}^{n} \frac{(-1)^k 4^{n-k} (2n-k)!}{k! (2n+1-2k)!} x^{n-k}.
$$
Note that the $n$ roots of $E_n(x)$ are  $\sin^2(\frac{j\,\pi}{2n+1})$, for $j=1, \dots,n$.
Define the {\it reciprocal polynomial} $\bar{E}_n(x):=x^nE_n(x^{-1})$. Then we have
$$
\bar{E}_n(x) \, = \, (-1)^n \, (2n+1) \sum_{k=0}^{n} \frac{(-1)^k 4^{k} (n+k)!}{(n-k)! (2k+1)!} x^{n-k}.
$$
Note that the $n$ roots of $\bar{E}_n(x)$ are  $\csc^2(\frac{j\,\pi}{2n+1})$.
It follows immediately that, denoting, as usual, the degree $k$ elementary symmetric function of
$\alpha_1, \dots, \alpha_n$ by $e_k(\alpha_1, \dots, \alpha_n)$, that

\bigskip

\noindent{\bf Lemma Top}
$$
e_k(\{\sin^2 (j \pi/(2n+1)), j=1 \dots n\}) \, = \,
\frac{4^{-k} \left(2 n -k \right)! \left(2 n +1\right)}{k ! \left(2 n -2 k +1\right)!}\cdot
$$

\bigskip

\noindent{\bf Lemma Ton}
$$
e_k(\{\csc^2(j \pi/(2n+1)), j=1 \dots n\}) \, = \,
\frac{\left(n +k \right)! 4^{k}}{\left(n -k \right)! \left(2 k +1\right)!}\cdot
$$

\medskip

Now we use Newton's identities \cite{newton} (see \cite{zeil-newton} 
for a lovely combinatorial proof). Let  
$$
p_k(\alpha_1, \dots, \alpha_n):= \sum_{j=1}^{n} \alpha_j^k
$$ 
be the {\it power-sum functions}, then
$$
k e_k (\alpha_1, \dots, \alpha_n) )\, = \,
\sum_{i=1}^k (-1)^{i-1} e_{k-i}(\alpha_1, \dots, \alpha_n)  p_i(\alpha_1, \dots, \alpha_n).
$$
This enables us (and our computers) to {\bf recursively} find explicit expressions for the power 
sums of both  $\{\sin^2(j \pi/(2n+1)), j=1 \dots n)\}$ and $\{\csc^2(j \pi/(2n+1)), j=1 \dots n)\}$.
Of course, for the former case we do not need it, since we have Theorem {\bf Top} above but it is 
still nice to know that it agrees up to $k=60$.

\bigskip

For types {\bf Tep} and {\bf Ten}, the underlying polynomial is $T_{2n}(\sqrt{x})$ whose roots are 
$$
\{\sin^2((2j+1)\pi/(4n)), j=1 \dots n\},
$$ 
and we have, analogously

\bigskip

\noindent{\bf Lemma Tep}
$$
e_k(\{\sin^2((2j+1)\,\pi/(4n)),j=1 \dots n\} ) \, = \,
\frac{2n\,4^{-k} \left(2 n -k -1\right)!}{k ! \left(2 n -2 k \right)!} \cdot
$$

\medskip

\noindent{\bf Lemma Ten}
$$
e_k(\{\csc^2 ((2j+1)\,\pi/(4n)),j=1 \dots n\} )\, = \,
\frac{n \left(n +k -1\right)! 4^{k}}{\left(n -k \right)! \left(2 k \right)!} \cdot
$$

\bigskip

For types {\bf Uop} and {\bf Uon}, the underlying polynomial is $U_{2n-1}(\sqrt{x})/\sqrt{x}$ 
whose roots are 
$$
\{\sin^2(j\pi/(2n)), j=1 \dots n-1\},
$$ 
and we have, analogously

\bigskip

\noindent{\bf Lemma Uop}
$$
e_k(\{\sin^2(j \pi/(2n)),j=1 \dots n-1\} )\, = \,
\frac{4^{-k} \left(2 n - k -1\right)!}{\left(2 n - 2 k -1\right)! k !} \cdot
$$

\medskip

\noindent{\bf Lemma Uon}
$$
e_k(\{\csc^2(j \pi/(2n)), j=1 \dots n-1\}) \, = \,
\frac{4^{k} \left(n +k \right)! \left(n -1\right)!}{\left(2 k +1\right)! \left(n -k -1\right)! n !}  \cdot
$$

\bigskip

For types {\bf Uep} and {\bf Uen}, the underlying polynomial is $U_{2n}(\sqrt{x})$ whose roots are 
$$
\{ \sin^2 ((2j-1)\pi/(4n+2)), j=1 \dots n\},
$$ 
and we have, analogously

\bigskip

\noindent{\bf Lemma Uep}
$$
e_k(\{\sin^2((2j-1)\,\pi/(4n+2)), j=1 \dots n\} )\, = \,
\frac{4^{-k} \left(2 n -k \right)!}{\left(2 n -2 k \right)! k !} \cdot
$$

\noindent{\bf Lemma Uen}
$$
e_k(\{\csc^2 ((2j-1) \pi/(4n+2)), j=1 \dots n\}) \, = \,
\frac{\left(n + k \right)! 4^{k}}{\left(2 k \right)! \left(n - k \right)!} \cdot
$$

\medskip

 From these, combined with Netwon's identities one can find as many sum-of-powers as one wishes.

\bigskip
\bigskip

\noindent
{\bf Acknowledgments.} \ \ We warmly thank D. Berend, M. Dekking, I. Gessel, J. Jorgenson,
B. Rand\'e and the three referees for their remarks and comments. In particular (added 
September 18, 2022), we would like to point to the following references that are related to our work:

\medskip

\small{

I. M. Gessel, Generating functions and generalized Dedekind sums, in
The Wilf Festschrift (Philadelphia, PA, 1996),
{\it Electron. J. Combin.} {\bf 4} (1997), Research Paper 11, 17 pp.

\medskip

C. A. Cadavid, P. Hoyos, J. Jorgenson, L. Smajlovi\'c, J. D. V\'elez,
On an approach for evaluating certain trigonometric character sums using 
the discrete time heat kernel, Preprint (2022), \url{https://arxiv.org/abs/2201.07878}.
}

\newpage

\begin{center}
------------------------
\end{center}

\begin{center}
{\Large Appendix by Shalosh B. Ekhad}
\end{center}

\begin{center}
------------------------
\end{center}

\noindent{\bf Type Ton}

\bigskip

\noindent{\bf Proposition Ton$_{\mathbf 1}$}
$$
\sum_{j=1}^{n} \csc^2\left(\frac{\pi j}{2n+1}\right) \, = \,
\frac{2 \left(n +1\right) n}{3} \cdot
$$

\noindent{\bf Proposition Ton$_{\mathbf 2}$}
$$
\sum_{j=1}^{n} \csc^4\left(\frac{\pi j}{2n+1}\right) \, = \,
\frac{8 \left(n +1\right) n \left(n^{2}+n +3\right)}{45}
$$

\noindent{\bf  Proposition Ton$_{\mathbf 3}$}
$$
\sum_{j=1}^{n} \csc^6\left(\frac{\pi j}{2n+1}\right) \, = \,
\frac{8 \left(n +1\right) n \left(8 n^{4}+16 n^{3}+35 n^{2}+27 n +54\right)}{945} \cdot
$$

\noindent{\bf  Proposition Ton$_{\mathbf 4}$}
$$
\sum_{j=1}^{n} \csc^8\left(\frac{\pi j}{2n+1}\right) \, = \,
\frac{128 \left(n +1\right) n \left(n^{2}+n +3\right) \left(3 n^{4}+6 n^{3}+7 n^{2}+4 n +15\right)}{14175} \cdot
$$

\noindent{\bf  Proposition Ton$_{\mathbf 5}$}
$$
\sum_{j=1}^{n} \csc^{10}\left(\frac{\pi j}{2n+1}\right) \, = \,
\frac{64 \left(n +1\right) n \left(16 n^{8}+64 n^{7}+182 n^{6}+322 n^{5}+493 n^{4}+524 n^{3}
+579 n^{2}+360 n +540\right)}{93555} \cdot
$$

For Propositions {\bf Ton$_{\mathbf k}$} for $6 \leq k \leq 50$, see the output file
$$
\text{\url{https://sites.math.rutgers.edu/~zeilberg/tokhniot/oTrigSums1.txt}}.
$$

\bigskip

\noindent{\bf Type Ten}

\bigskip

\noindent{\bf  Proposition Ten$_{\mathbf 1}$}
$$
\sum_{j=1}^{n} \csc^2\left(\frac{\pi (2j-1)}{4n}\right) \, = \,
2 n^2 .
$$

\noindent{\bf  Proposition Ten$_{\mathbf 2}$}
$$
\sum_{j=1}^{n} \csc^4\left(\frac{\pi (2j-1)}{4n}\right) \, = \,
\frac{4 n^{2} \left(2 n^{2}+1\right)}{3} \cdot
$$

\noindent{\bf  Proposition Ten$_{\mathbf 3}$}
$$
\sum_{j=1}^{n} \csc^6\left(\frac{\pi (2j-1)}{4n}\right) \, = \,
\frac{8 n^{2} \left(8 n^{4}+5 n^{2}+2\right)}{15} \cdot
$$

\noindent{\bf  Proposition Ten$_{\mathbf 4}$}
$$
\sum_{j=1}^{n} \csc^8\left(\frac{\pi (2j-1)}{4n}\right) \, = \,
\frac{16 n^{2} \left(136 n^{6}+112 n^{4}+49 n^{2}+18\right)}{315} \cdot
$$

\noindent{\bf  Proposition Ten$_{\mathbf 5}$}
$$
\sum_{j=1}^{n} \csc^{10}\left(\frac{\pi (2j-1)}{4n}\right) \, = \,
\frac{32 n^{2} \left(992 n^{8}+1020 n^{6}+546 n^{4}+205 n^{2}+72\right)}{2835} \cdot
$$

\bigskip

For Propositions {\bf Ten$_{\mathbf k}$} for $6 \leq k \leq 50$, see the output file
$$
\text{\url{https://sites.math.rutgers.edu/~zeilberg/tokhniot/oTrigSums5.txt}}.
$$

\bigskip

\noindent{\bf Type Uon}

\bigskip

\noindent{\bf  Proposition Uon$_{\mathbf 1}$}
$$
\sum_{j=1}^{n-1} \csc^2\left(\frac{j \pi }{2n}\right) \, = \,
\frac{2 \left(n +1\right) \left(n -1\right)}{3} \cdot
$$

\noindent{\bf  Proposition Uon$_{\mathbf 2}$}
$$
\sum_{j=1}^{n-1} \csc^4\left(\frac{j \pi }{2n}\right) \, = \,
\frac{4 \left(n +1\right) \left(n -1\right) \left(2 n^{2}+7\right)}{45} \cdot
$$

\noindent{\bf  Proposition Uon$_{\mathbf 3}$}
$$
\sum_{j=1}^{n-1} \csc^6\left(\frac{j \pi }{2n}\right) \, = \,
\frac{8 \left(n +1\right) \left(n -1\right) \left(8 n^{4}+29 n^{2}+71\right)}{945} \cdot
$$

\noindent{\bf  Proposition Uon$_{\mathbf 4}$}
$$
\sum_{j=1}^{n-1} \csc^8\left(\frac{j \pi }{2n}\right) \, = \,
\frac{16 \left(n +1\right) \left(n -1\right) \left(24 n^{6}+104 n^{4}+251 n^{2}+521\right)}{14175} \cdot
$$

\noindent{\bf  Proposition Uon$_{\mathbf 5}$}
$$
\sum_{j=1}^{n-1} \csc^{10}\left(\frac{j \pi }{2n}\right) \, = \,
\frac{32 \left(n +1\right) \left(n -1\right) \left(32 n^{8}+164 n^{6}+450 n^{4}+901 n^{2}+1693\right)}{93555}
\cdot
$$

\bigskip

For Propositions {\bf Uon$_{\mathbf k}$} for $6 \leq k \leq 50$, see the output file
$$
\text{\url{https://sites.math.rutgers.edu/~zeilberg/tokhniot/oTrigSums9.txt}}.
$$

\bigskip

\noindent{\bf Type Uen}

\bigskip

\noindent{\bf  Proposition Uen$_{\mathbf 1}$}
$$
\sum_{j=1}^{n} \csc^2\left(\frac{(2j-1) \pi }{4n+2}\right) \, = \,
2 \left(n +1\right) n .
$$

\noindent{\bf  Proposition Uen$_{\mathbf 2}$}
$$
\sum_{j=1}^{n} \csc^4\left(\frac{(2j-1) \pi }{4n+2}\right) \, = \,
\frac{8 \left(n +1\right) n \left(n^{2}+n +1\right)}{3} \cdot
$$

\noindent{\bf  Proposition Uen$_{\mathbf 3}$}
$$
\sum_{j=1}^{n} \csc^6\left(\frac{(2j-1) \pi }{4n+2}\right) \, = \,
\frac{8 \left(n +1\right) n \left(8 n^{4}+16 n^{3}+19 n^{2}+11 n +6\right)}{15} \cdot
$$

\noindent{\bf  Proposition Uen$_{\mathbf 4}$}
$$
\sum_{j=1}^{n} \csc^8\left(\frac{(2j-1) \pi }{4n+2}\right) \, = \,
\frac{128 \left(n +1\right) n \left(n^{2}+n +1\right) \left(17 n^{4}+34 n^{3}+31 n^{2}+14 n +9\right)}{315}
\cdot
$$

\noindent
{\bf  Proposition Uen$_{\mathbf 5}$}
$$
\sum_{j=1}^{n} \csc^{10}\left(\frac{(2j-1) \pi }{4n+2}\right) \, 
$$
$$
= \,
\frac{64 \left(n +1\right) n \left(496 n^{8}+1984 n^{7}+4106 n^{6}+5374 n^{5}+4979 n^{4}
+3316 n^{3}+1669 n^{2}+576 n +180\right)}{2835} \cdot
$$

\bigskip

For Propositions {\bf Uen$_{\mathbf k}$} for $6 \leq k \leq 50$, see the output file
$$
\text{\url{https://sites.math.rutgers.edu/~zeilberg/tokhniot/oTrigSums12.txt}}.
$$

\end{document}